\documentclass[11pt,leqno]{article} 
\usepackage{graphics}
\newtheorem{thm}{Theorem}[section]
\newtheorem{lma}{Lemma}[section]

\newcommand{\beqa}{\begin{eqnarray}}
\newcommand{\eeqa}{\end{eqnarray}}

\newcommand{\pf}{\noindent {\bf Proof:} $\s$ }
\newcommand{\epf}{ \hfill$\diamondsuit$ \medskip}

\newcommand{\ds}{\displaystyle}
\newcommand{\beq}{\begin{equation}}
\newcommand{\eeq}{\end{equation}}
\newcommand{\lbl}{\label}
\newcommand{\s}{\; \;}

\newcommand{\ep}{\epsilon}

\newcommand{\ra}{\rightarrow}
\newcommand{\al}{\alpha}
\newcommand{\p}{\varphi}
\title{Regularization of radial solutions of $p$-Laplace equations, and computations using infinite series  }

\author{
Philip Korman   \\ 
Department of Mathematical Sciences \\ 
University of Cincinnati \\ 
Cincinnati Ohio 45221-0025 \\
kormanp@ucmail.uc.edu\\
}

\date{}

\begin{document}

\maketitle
\begin{abstract} 
We consider radial solutions of equations with the $p$-Laplace operator in $R^n$. We introduce a change of variables, which in effect removes the singularity at $r=0$. While solutions are not of class $C^2$, in general, we show that solutions are $C^2$ functions of $\ds r^{\frac{p}{2(p-1)}}$. Then we express the solution as an infinite series in powers of $\ds r^{\frac{p}{p-1}}$, and give explicit formulas for its coefficients. We implement this algorithm, using {\em Mathematica}. {\em Mathematica}'s ability to perform the exact computations turns out to be crucial.
 \end{abstract}

\begin{flushleft}
\end{flushleft}

\begin{flushleft}
\end{flushleft}

\section{Introduction}
\setcounter{equation}{0}
\setcounter{thm}{0}
\setcounter{lma}{0}
Recently there has been an enormous interest in  equations with the $p$-Laplace operator  in $R^n$ (with $p>1$, $u=u(x)$, $ x \in R^n$)
\[
\mbox{div} \left( |\nabla u|^{p-2} \nabla u \right) +f(u)=0 \,,
\]
see e.g., a review by P. Drabek \cite{D}, and the important paper of B. Franchi et al \cite{FLS}.
Radial solutions of this equation, with the initial data at $r=0$, satisfy
\beq
\lbl{1}
\p (u')'+\frac{n-1}{r} \, \p (u') +f(u)=0, \s u(0)=\al >0, \; u'(0)=0 \,,
\eeq
where $\p (v)=v |v|^{p-2}$, $p>1$. To guess the form of the solution, let us drop the higher order term and consider
\beq
\lbl{2}
\frac{n-1}{r} \, \p (u') +f(u)=0, \s u(0)=\al  \,.
\eeq
This is a completely different equation, however, in case $p=2$, it is easy to check that the  {\em form} of  solutions is the same: in both cases, it is a series $\sum _{n=0}^{\infty} a_n r^{2n}$ (with different coefficients), see P. Korman \cite{K4} or \cite{K1}. It is natural to guess that the form of solutions will be same for (\ref{1}) and (\ref{2}), in case $p \ne 2$ too. With that in mind, let us solve (\ref{2}) in case $f(u)=e^u$. Since $\al >0$, we see from (\ref{2}) that  $u'(r)<0$ for  all $r$. Then $\p (u')=-(-u')^{p-1}$, and we have 
\[
u'=-\left[\frac{r}{n-1} e^u \right]^{\frac{1}{p-1}}=- \frac{r^{\frac{1}{p-1}}}{(n-1)^{\frac{1}{p-1}}}e^{\frac{1}{p-1}u} \,.
\]
Integrating, we get
\[
u(r)=\al-(p-1) \ln \left(1+\frac{e^{\frac{1}{p-1} \al}}{p (n-1)^{\frac{1}{p-1}}}\, r^{\frac{p}{p-1}}\right) \,.
\]
We see that $u(r)$ is a function of $r^{\frac{p}{p-1}}$, and, for $r $ small, we can expand it as a series $u(r)=\sum_{n=0}^{\infty} b_n r^{n \frac{p}{p-1}}$, with some coefficients $b_n$. Motivated by this example, we make a change of variables $r \ra z$ in (\ref{1}), by letting $z^2=r^{\frac{p}{p-1}}$. We expect solutions of (\ref{1}) to be of the form $\sum_{n=0}^{\infty} c_n z^{2n}$, which is a real analytic function of $z$, if this series converges. \medskip

The following lemma provides the crucial change of variables.
\begin{lma}\lbl{lma:1}
Denote 
\beq
\lbl{2a}
\s\s \bar \al =\frac{p}{2(p-1)}, \s  \beta=\frac{1}{\bar \al} (\bar \al -1)=\frac{-p+2}{p}, \s  \gamma=\beta(p-1)=3-p-\frac{2}{p} \,,
\eeq
\[
 a=\bar \al ^p, \s   A=\bar \al ^p \gamma +(n-1)\bar \al ^{p-1} \,.
\] 
Then, for $p>2$, the change of variables $z^2=r^{\frac{p}{p-1}}$ transforms (\ref{1}) into
\beq
\lbl{3}
 au''(z)+\frac{A}{(p-1)z}u'(z)+\frac{z^{p-2}}{\p '(u'(z))} f(u)=0, \s u(0)=\al,   \; u'(0)=0 \,.
\eeq 
Conversely, if the solution of (\ref{3}) is of the form $u=v(z^2)$, with $v(t) \in C^1(\bar R_+)$, then the same change of variables transforms (\ref{3}) into (\ref{1}), for any $p>1$.
\end{lma}

\pf
We have $\ds z=r^{\al}$, $\ds  \frac{du}{dr}=\al \frac{du}{dz}r^{\al -1} =\al z^{\beta} \frac{du}{dz} $. By the homogeneity  of $\p$, we have ($\p (cv)=c^{p-1} \p (v)$, for any $c>0$)
\[
\p (u'(r)) =\p \left( \al z^{\beta} u_z \right)=\al ^{p-1} z^{\beta (p-1)} \p (u_z)=\al ^{p-1} z^{\gamma} \p (u_z) \,.
\]
Then (\ref{1}) becomes
\[
\al z^{\beta} \frac{d}{dz} \left[\al ^{p-1} z^{\gamma} \p (u_z) \right]+(n-1)z^{-2+2/p} \al ^{p-1} z^{\gamma} \p (u_z)+f(u)=0 \,,
\]
which simplifies to
\[
az \p'(u'(z)) u''(z)+A \p(u'(z))+z^{1/\al -\gamma} f(u)=0 \,.
\]
This implies (\ref{3}), keeping in mind that $\p'(v)=(p-1)|v|^{p-2}$, $\p(v)=\frac{1}{p-1}v \p'(v)$, and that $1/\al -\gamma-1=p-2$.
Also, $\ds \frac{du}{dz}=\frac{2(p-1)}{p}\frac{du}{dr} \, z^{\frac{p-2}{p}}$, so that $\frac{du}{dz}(0)=0$. (It is only on the last step that we need $p>2$.)
\medskip

Conversely, our change of variables transforms the equation in (\ref{3}) into the one in (\ref{1}). Under our assumption, $u(r)=v(r^{\frac{p}{p-1}})$, so that $u'(0)=0$ for any $p>1$.
\epf

The change of variables $z^2=r^{\frac{p}{p-1}}$ in effect removes the singularity at zero for $p$-Laplace equations. Indeed, 
\[
 \lim _{z \ra 0}  \frac{z^{p-2}}{\p '(u'(z))}=\frac{1}{(p-1) |u''(0)|^{p-2}} \,, 
\]
which lets us compute $u''(0)$ from the equation (\ref{3}) (the existence of $u''(0)$ is proved later). Indeed, assuming that $f(\al)>0$, we have 
\beq
\lbl{3.1}
u''(0)=-\left[ \frac{f(\al)}{a(p-1)+A} \right]^{\frac{1}{p-1}} \,.
\eeq
(In case $f(\al)<0$, we have $u''(0)=\left[ \frac{-f(\al)}{a(p-1)+A} \right]^{\frac{1}{p-1}}$.)
 We prove that  $u(z)$ is smooth, provided that $f(u)$ is smooth. It follows that the solution of $p$-Laplace problem (\ref{1})  has the form $u(r^{\frac{p}{2(p-1)}})$, with smooth $u(z)$. We believe that our reduction of  the $p$-Laplace equation (\ref{2})  to the form (\ref{3}) is likely to find other applications.
\medskip

We express the solution of (\ref{1}) in the form $u(r)=\sum_{k=0}^{\infty} a_k r^{k \frac{p}{p-1}}$, and present {\em explicit } formulas to compute the coefficients $a_k$. Interestingly, the coefficient $a_1$ turned out to be special, as it enters in two ways the formula for other $a_k$. Our formulas are easy to implement in {\em Mathematica}, and very accurate series approximations can be computed reasonably quickly. We utilize {\em Mathematica}'s ability to perform the ``exact computations", as we explain in Section $3$.

\section{Regularity of solutions in case $p>2$}
\setcounter{equation}{0}
\setcounter{thm}{0}
\setcounter{lma}{0}

It is well known that solutions of $p$-Laplace equations are not of class $C^2$, in general. In fact, rewriting the equation in (\ref{1}) in the form
\beq
\lbl{+}
(p-1)u''+\frac{n-1}{r}u'+|u'|^{2-p}f(u) =0\,,
\eeq
we see that in case $p>2$, $u''(0)$ does not exist. We show that in this case the solution of (\ref{1}) is a $C^2$ function of $r^{\frac{p}{2(p-1)}}$.

We rewrite the equation in  (\ref{1})  as
\beq
\lbl{*}
r^{n-1} \p (u'(r))=-\int_0^r t^{n-1} f(u(t)) \, dt \,.
\eeq
Observe that $ \p^{-1}(t)=-(-t)^{^{\frac{1}{p-1}}}$, for $t<0$. If we assume that $f(\al)>0$, then for small $r>0$, we may express from (\ref{*})
\beq
\lbl{***}
-u'(r)=\frac{1}{r^{\frac{n-1}{p-1}}} \left[ \int _0^r t^{n-1} f(u(t)) \, dt \right]^{^{\frac{1}{p-1}}} \,.
\eeq
Integrating
\beq
\lbl{**}
u(r)=\al - \int _0^r \frac{1}{t^{\frac{n-1}{p-1}}} \left[ \int _0^t s^{n-1} f(u(s)) \, ds \right]^{^{\frac{1}{p-1}}} \, dt\,.
\eeq

We recall the following lemma from J.A. Iaia \cite{I}.
\begin{lma}\lbl{lma:0}
Assume that $f(u)$ is Lipschitz continuous. Then one can find an $\epsilon >0$, so that the problem (\ref{1}) has a unique solution $u(r) \in C^1[0,\epsilon)$. In case $1<p \leq 2$, $u(r) \in C^2[0,\epsilon)$. 
\end{lma}

In the space $C[0,\epsilon)$ we denote $B_R^{\ep}=\{u \in  C[0,\epsilon)$, such that $||u-\al|| \leq R \}$, where $||\cdot||$ denotes the norm in $C[0,\epsilon)$. 
The proof of Lemma \ref{lma:0} involved showing  that the map $T(u)$, defined by the right hand side of (\ref{**}), is a contraction, taking $B_R^{\ep}$ into itself, for any $R>0$, and $\ep$ sufficiently small (see \cite{I}, and also  \cite{PS1} for a similar argument). This argument provided a continuous solution of (\ref{**}), which by  (\ref{***}) is in  $C^1[0,\epsilon)$, and in case  $1<p \leq 2$, $u(r) \in C^2[0,\epsilon)$, by (\ref{+}) (from (\ref{***}) it follows that the limit $\lim _{r \ra 0} \frac{u'(r)}{r}=u''(0)=0$ exists).
\epf

In case  $p>2$,  we have the following regularity result.
\begin{thm}\lbl{thm:2}
Assume that $p>2$, $f(u)$ is Lipschitz continuous and $f(\al)>0$. For $\ep >0$ sufficiently small, the problem (\ref{1}) has a  solution of the form $u \left( r^{\frac{p}{2(p-1)}} \right)$, where $u(z) \in C^2[0,\epsilon)$, and $u'(z)<0$ on $(0,\ep)$. This solution is unique among all continuous functions satisfying (\ref{**}).
If, moreover, $f(u) \in C^k$, then $u(z) \in C^{k+2}[0,\epsilon)$.
\end{thm}

\pf
By Lemma \ref{lma:0} we have a unique solution of the problem (\ref{1}), $u(r) \in C^1[0,\epsilon _1)$, for some $\epsilon _1>0$ small.
By Lemma \ref{lma:1}, this translates to a solution of the problem (\ref{3}), $u(z) \in C^1[0,\epsilon _1)$.
With $m=\frac{A}{a(p-1)}$, we multiply the equation in  (\ref{3}) by $z^m$, and rewrite it as 
\[
-\frac{u'(z)}{z}=\frac{1}{a(p-1)} \frac{\int _0^z \frac{t^{m+p-2}}{|u'(t)|^{p-2}} f(u(t)) \,dt}{z^{m+1}} \,.
\]
Taking the limit as $z \ra 0$, and denoting $L=\lim _{ z \ra 0} \frac{u'(z)}{z}$, we get
\[
-L=\frac{f(\al)}{a(p-1)(m+1)|L|^{p-2}} \,.
\]
It follows that this limit $L$ exists, proving the existence of  $u''(0)$, as given by (\ref{3.1}). Observe that $u''(0)<0$. It follows that $u'(z)<0$ for small $z$, so that $\p ' (u'(z))<0$, and then $u(z) \in C^2[0,\epsilon)$, from the equation (\ref{3}).
\medskip

Assume that $f(u) \in C^1$. Differentiate the equation (\ref{3})
\[
au'''+\frac{A}{p-1}\frac{u'' z-u'}{z^2}+\frac{p-2}{p-1} \left(-\frac{z}{u'} \right)^{p-3} \frac{zu''-u'}{{u'}^2} f(u)+\frac{1}{p-1} \left(-\frac{z}{u'} \right)^{p-2}f' u'=0 \,.
\]
From here, $u(z) \in C^3(0,\epsilon)$. Letting $z \ra 0$, and using that $\lim _{z \ra 0} \frac{u'' z-u'}{z^2}=\frac12 u'''(0)$, and $\lim _{z \ra 0} \frac{zu''-u'}{{u'}^2}=\frac{u'''(0)}{2{u''(0)}^2}$, we conclude that $u'''(0)=0$ (the existence of $u'''(0)$ is proved as before). It follows that $u(z) \in C^3[0,\epsilon)$. Higher regularity is proved by taking further derivatives of the equation.
\epf

\section{Representation of solutions using infinite series}
\setcounter{equation}{0}
\setcounter{thm}{0}
\setcounter{lma}{0}

We shall consider an auxiliary problem
\beq
\lbl{5}
 au''(z)+\frac{A}{(p-1)z}u'(z)+\frac{|z|^{p-2}}{\p '(u'(z))} f(u)=0, \s u(0)=\al,   \; u'(0)=0 \,.
\eeq

\begin{lma}\lbl{lma:3}
Any solution of the problem (\ref{5}) is an even function.
\end{lma}

\pf
Observe that the change of variables $z \ra -z$ leaves (\ref{5}) invariant. If solution $u(z)$ were not even, then $u(-z)$ would be another solution of (\ref{5}). By Lemma \ref{lma:1}, $u(z)$ and $u(-z)$ translate into two different solutions of the problem (\ref{1}), contradicting the uniqueness part of Lemma \ref{lma:0}.
\epf

It follows from the last lemma that any series solution of (\ref{5}) must be of the form $\sum _{n=0}^{\infty} a_n z^{2n}$. The same must be true for the problem (\ref{3}), since for $z>0$ it agrees with (\ref{5}). Numerically, we shall be computing the partial sums
$\sum _{n=0}^{k} a_n z^{2n}$, which will provide us with the solution, up to the terms of order $O(z^{2k+2})$. Write the partial sum in the form
\beq
\lbl{6}
u(z)=\bar u(z)+a _{k} z^{2k} \,,
\eeq
where $\bar u(z)=\sum _{n=0}^{k-1} a_n z^{2n}$. We regard $\bar u(z)$ as already computed, and the question is how to compute $a_k$.
Using the  constants defined in (\ref{2a}), we let
\[
B_k= 2k(2k-1)a+\frac{2kA}{p-1}  \,.
\]
\begin{thm}\lbl{thm:1}
Assume that $\al >0$, $f(u) \in C^{\infty}(R)$, and $f(\al)>0$. The solution of the problem (\ref{3}) in terms of a series of the form $\sum _{k=0}^{\infty} a_k z^{2k}$ is obtained by taking $a_0=\al$, then 
\beq
\lbl{9}
a_1=- \left[ \frac{1}{(p-1)2^{p-2}B_1} f(\al) \right]^{\frac{1}{p-1}} <0 \,,
\eeq
and for $k \geq 2$, we have  (the following limits exist)
\beq
\lbl{10}
a_k=-\frac{1}{B_k C_k} \lim _{z \ra 0} \frac{\frac{z^{p-2}}{\p '(\bar u'(z))} f(\bar u)+a \bar u''(z)+\frac{A}{(p-1)z}\bar u'(z)}{z^{2k-2}} \,,
\eeq
where $\bar u=\sum _{n=0}^{k-1} a_n z^{2n}$ is the previously computed approximation, and 
\[
C_k=1+\frac{ k(p-2) f(\al )}{(p-1)2^{p-2} B_k (-a_1)^{p-1}} \,.
\]
\end{thm}

\pf
Plugging $u=\al+a_1 z^2$ into the equation (\ref{3}), gives
\[
a_1B_1=-\frac{z^{p-2}}{(p-1) |2a_1z|^{p-2}} f(\al+a_1 z^2) \,.
\]
Letting $z \ra 0$, 
\[
a_1 B_1 =-\frac{1}{(p-1) |2a_1|^{p-2}} f(\al) \,,
\]
which implies that $a_1<0$, leading to (\ref{9}). Of course, $u=\al+a_1 z^2$ is not a solution of (\ref{3}). But the other terms of the solution $u(z)=\sum _{k=0}^{\infty} a_k z^{2k}$ produce a correction, which disappears in the limit. Indeed, plugging $u(z)=\al+a_1z^2+\sum _{k=2}^{\infty} a_k z^{2k}$ into  (\ref{3}), gives
\[
a_1B_1+\sum _{k=2}^{\infty} a_k B_k z^{2k-2} 
\]
\[
=-\frac{z^{p-2}}{(p-1) |2a_1z+\sum _{k=2}^{\infty} 2ka_k z^{2k-1}|^{p-2}} f ( \al+a_1 z^2+\sum _{k=2}^{\infty} a_k z^{2k} ) \,,
\]
and going to the limit, with $z \ra 0$, gives the same value of $a_1$.
\medskip

Plugging  $u(z)=\bar u(z)+a _{k} z^{2k}$ into the equation (\ref{3}), gives
\beq
\lbl{12}
-a_k B_k =\frac{\frac{z^{p-2}}{\p'(\bar u'+2ka_kz^{2k-1})}f( \bar u(z)+a _{k} z^{2k}) +a \bar u''(z)+\frac{A}{(p-1)z} \, \bar u'(z)}{z^{2k-2}} \,.
\eeq
We now expand the quotient in the first term in the numerator. In this expansion we do not need to show the terms that are of order $O(z^{2k-1})$ and higher, since 
$\lim _{z \ra 0} \frac{O(z^{2k-1})}{z^{2k-2}}=0$. For $z>0$ and small, we have (observe that $u'(z)<0$, and $u'(z) \sim 2a_1 z$, for $z$ small)
\[
(p-1)\frac{z^{p-2}}{\p'(\bar u'+2ka_kz^{2k-1})} =\frac{z^{p-2}}{\left(-\bar u'-2ka_kz^{2k-1}\right)^{p-2}}=
\]
\[
\frac{z^{p-2}}{\left(-\bar u' \right)^{p-2} \left(1+2ka_k \frac{z^{2k-1}}{\bar u'(z)}\right)^{p-2}}=\frac{z^{p-2}}{\left(-\bar u' \right)^{p-2}}
\left(1-2ka_k(p-2) \frac{z^{2k-1}}{\bar u'(z)} \right)+O(z^{2k-1})
\]
\[
=\frac{z^{p-2}}{\left(-\bar u' \right)^{p-2}} \left(1+2ka_k(p-2) \frac{z^{2k-2}}{(-2a_1)} \right)+O(z^{2k-1}) \,
\]
\[
=\frac{(p-1)z^{p-2}}{\p'(\bar u' )}+2ka_k(p-2) \frac{z^{2k-2}}{(-2a_1)(-2a_1)^{p-2}}+O(z^{2k-1}) \,.
\]
(Observe that $\ds \frac{z}{- \bar u'}=\frac{1}{-2 a_1}+o(z)$, and  $\ds \frac{z^{p-2}}{(- \bar u')^{p-2}}=\frac{1}{(-2a_1)^{p-2}}+o(z)$.)
Also
\[
f(\bar u(z)+a_k z^{2k})=f(\bar u(z))+O(z^{2k}) \,.
\]
Using these expressions in (\ref{12}), and taking the limit, we get
\[
a_k=-\frac{1}{B_k } \lim _{z \ra 0} \frac{\frac{z^{p-2}}{\p '(\bar u'(z))} f(\bar u)+a \bar u''(z)+\frac{A}{(p-1)z}\bar u'(z)}{z^{2k-2}} -\frac{k(p-2) f(\al)}{(p-1)B_k2^{p-2} (-a_1)^{p-1}} \, a_k\,.
\]
(By Theorem \ref{thm:2}, $u(z) \in C^{\infty}[0,\ep)$. Hence, the limit representing $a_k=\frac{u^{2k}(0)}{(2k)!}$ exists.)
Solving this equation for $a_k$, we conclude (\ref{10}). Plugging $u(z)=\bar u(z)+a_k z^{2k}+\sum _{n=k+1}^{\infty} a_nz^{2n}$ into (\ref{3}), produces the same formula for $a_k$.
\epf

With $a_k$'s computed as in this theorem, the series $\sum _{k=0}^{\infty} a_{k} r^{\frac{kp}{p-1}}$ gives the solution to the original problem (\ref{1}). In case $p=2$, we proved in \cite{K4} that when $f(u)$ is real analytic, the series $\sum _{k=0}^{\infty} a_{k} z^{2k}$ converges for small $z$, giving us a real analytic solution. It is natural to expect convergence for $p \ne 2$ too, so that the solution of (\ref{1}) is a real analytic function of $r^{\frac{p}{p-1}}$.

\section{Numerical computations}
\setcounter{equation}{0}
\setcounter{thm}{0}
\setcounter{lma}{0}

It is easy to implement our formulas for computing the solution in {\em Mathematica}. It is crucial that {\em Mathematica} can perform exact computations for fractions. If one tries floating point computations, the limits in (\ref{10}) become infinite. All numbers must be entered as fractions. For example, one cannot enter $p=4.1$, it should be $p=\frac{41}{10}$ instead. ({\em Mathematica} switches to floating point computations, once it sees a  number  entered as a floating point.)
\medskip

\noindent
{\bf Example} We  solved
\beq
\lbl{14}
\p (u')'+\frac{n-1}{r} \, \p (u') +e^u=0, \s u(0)=1, \; u'(0)=0 \,,
\eeq
with $\p (v)=v |v|^{p-2}$, and $p=\frac{41}{10}$. {\em Mathematica}  calculated that the corresponding equation (\ref{3}) is
\beq
\lbl{15}
\s\s\s a(p-1) u''(z)+\frac{A}{z}u'(z)+\frac{z^{p-2}}{(-u'(z))^{p-2}} e^{u(z)}=0, \s u(0)=1,   \; u'(0)=0 \,,
\eeq
with $a(p-1)=\frac{2825761
   \sqrt[10]{\frac{41}{62}}}{476
   6560}$, $A=\frac{1309499
   \sqrt[10]{\frac{41}{62}}}{476
   6560}$.
When we computed the series solution of (\ref{15}) up to $a_5$, {\em Mathematica} returned (instantaneously)
\[
u(z)=1-\frac{31}{41}
   \left(\frac{e}{3}\right)^{10/
   31} z^2+\frac{4805 \,
   3^{11/31} \, e^{20/31}
   }{225254}\, z^4-\frac{326241241
   \left(\frac{e}{3}\right)^{30/
   31}
   }{43314091660}z^6
\]
\[
+\frac{513127652305
   79\,  e^{40/31}
   }{154203017487865920 \,
   3^{9/31}}\, z^8-\frac{13334484822273130589\,
   e^{50/31}
   }{283500239799651287332
   500 \,
   3^{19/31}} \, z^{10} \,.
\]
\medskip

\noindent
The same solution using floating point numbers is
\[
u(z)=1-0.732424 z^2+0.0600499
   z^4-0.00684643 z^6
\]
\[
+0.000879009
   z^8-0.000120356 z^{10} \,.
\]
For the original equation (\ref{14}), this implies (we have $\frac{p}{p-1}=\frac{41}{31} $, and $z^2=r^{\frac{41}{31}}$)
\[
u(r)=1-0.732424 r^{\frac{41}{31}}+0.0600499
   r^{\frac{82}{31}}-0.00684643 r^{\frac{123}{31}}
\]
\[
+0.000879009
  r^{\frac{164}{31}}-0.000120356 r^{\frac{205}{31}}+ \cdots \,.
\]
To check the accuracy of this computation, we  denoted by $q(z)$ the left hand side of (\ref{15}) (with $u(z)$ being the above polynomial of  tenth degree),
and asked {\em Mathematica} to expand $q(z)$ into series about $z=0$. {\em Mathematica} returned: $q(z)=O(z^{\frac{121}{10}})$. We have performed similar computations, with similar results, for other values of $u(0)=\al$. For larger values of $\al$, e.g., for  $\al=2$, the computations take longer, but no more than several minutes.
\medskip

We have obtained similar results for all other $f(u)$ and $p$ that we tried (including the case $1<p<2$). We wish to stress that in all computations, when the solution of (\ref{3}) was computed up to the order $z^{2n}$, the defect function $q(z)$ was at least of order $O(z^{2n+2})$ near $z=0$. This heuristic result is consistent with the Theorem \ref{thm:1}, but does not seem to follow from it.


\begin{thebibliography}{99}
\bibitem{D}
P. Drabek,  The $p$-Laplacian - mascot of nonlinear analysis, {\em   Acta Math. Univ. Comenian. (N.S.)}  {\bf 76},  no. 1, 85-98   (2007).
\vspace{-0.2cm}

\bibitem{FLS}
B. Franchi, E. Lanconelli  and J.  Serrin,  Existence and uniqueness of nonnegative solutions of quasilinear equations in $ R\sp n$, {\em Adv. Math.} {\bf 118}, no. 2, 177-243  (1996).
\vspace{-0.2cm}

\bibitem{I}
J.A. Iaia,  Localized solutions of elliptic equations: loitering at the hilltop, {\em  Electron. J. Qual. Theory Differ. Equ.}, No. 12, 15 pp. (electronic) (2006).
\vspace{-0.2cm}

\bibitem{K4}
P. Korman,  Computation of radial solutions of semilinear equations, {\em   Electron. J. Qual. Theory Differ. Equ.}, No. 13, 14 pp. (electronic) (2007).
\vspace{-0.2cm}

\bibitem{K1}
P. Korman, Global Solution Curves for Semilinear Elliptic Equations, World Scientific, Hackensack, NJ (2012).
 \vspace{-0.2cm}

\bibitem{PS1}
L.A. Peletier and J. Serrin, Uniqueness of positive solutions of semilinear equations in $ R\sp{n}$, {\em Arch. Rational Mech. Anal.} {\bf 81}  no. 2, 181-197 (1983).
\vspace{-0.2cm}

\end{thebibliography}
\end{document}